\documentclass[11pt]{article} 

\usepackage[utf8]{inputenc} 

\usepackage{amsmath}
\usepackage{amssymb}
\usepackage{amsthm}

\usepackage[pdftex]{graphicx}

\usepackage{geometry} 
\usepackage{graphicx} 

\usepackage{booktabs} 
\usepackage{array} 
\usepackage{paralist} 
\usepackage{verbatim} 
\usepackage{subfig} 

\usepackage{fancyhdr} 
\usepackage{sectsty}
\allsectionsfont{\sffamily\mdseries\upshape} 

\usepackage[nottoc,notlof,notlot]{tocbibind} 
\usepackage[titles,subfigure]{tocloft} 

\title{Distribution of Euclidean Distances Between \\* Randomly Distributed Gaussian Points in n-Space}
\author{Benjamin Thirey, M.S. and Randal Hickman, PhD \\*
Department of Mathematical Sciences \\*
United States Military Academy, West Point}
\date{} 

\begin{document}
\maketitle

\begin{abstract}
The curse of dimensionality is a common phenomenon which affects analysis of datasets characterized by large numbers of variables associated with each point.  Problematic scenarios of this type frequently arise in classification algorithms which are heavily dependent upon distances between points, such as nearest-neighbor and $k$-means clustering.  Given that contributing variables follow Gaussian distributions, this research derives the probability distribution that describes the distances between randomly generated points in n-space. The theoretical results are extended to examine additional properties of the distribution as the dimension becomes arbitrarily large.  With this distribution of distances between randomly generated points in arbitrarily large dimensions, one can then determine the significance of distance measurements between any collection of individual points. 
\end{abstract}

\section{Introduction}

    In recent years, there has been a significant increase in size among large datasets which possess hundreds or even thousands of features per observation.  This has given rise to increased observations of the phenomenon known as the "curse of dimensionality." First identified by Bellman in 1957, an application of this phenomenon occurs as the distances between observations grow with the number of additional dimensions. [3]  One of the side effects of this phenomenon is that spurious results are obtained due to the increase in variance produced by the introduction of extraneous features (dimensions).  These spurious results arise directly from the influence that additional distances have on the underlying metric used to calculate the nearest distance for classification examples in a training dataset.  This effect is especially problematic in statistical learning algorithms which rely heavily on measures of closeness or similarity between observations. \\
\indent
    Interest in the distances between randomly generated points was initially undertaken by trying to analyze the effects that additional data features (dimensions) have upon classification algorithms used for statistical learning with large datsets.  Hastie, et al. provide an in-depth analysis of this phenomenon in their text ``The Elements of Statistical Learning", which is a widely used textbook at the graduate level. [7]  The Nearest Neighbor algorithm is one such common machine learning technique which is known to be sensitive to the effects of extraneous dimensions because additional features increase the distances between points and reduce the contrast between meaningful neighbors and those which may arise purely from statistical noise.  Aggarwal, Hinneburg, and Keim provide a very good discussion of this phenomena in their paper and term it ``relative contrast". [1]  They observe that as the number of dimensions increases, the distance to a the nearest query point using a nearest neighbor algorithm increases faster than the difference in distances between the furthest and nearest query points.  Their conclusion is that as the number of dimensions increases, proximity queries become meaningless due to proper discrimination regarding what constitutes a suitable neighbor for classification purposes. \\
\indent    
While significant research has considered this problem conceptually, the authors are aware of little research regarding a quantitative approach.  Several years ago, there were several successful attempts to determine the distribution of distances between randomly generated points uniformly distributed within a hypersphere [2, 6, 9], randomly generated points distributed on the surface of a hypersphere [10], or more general analysis of random points contained within some compact convex subset of points in Euclidean space [4].  Solomon provided additional extensive coverage of these concepts for uniformly distributed points in circles and spheres in his book titled ``Geometric Probability." [11]  Interestingly, all of these observations use a uniform distribution of point density, either on the surface or interior of a sphere and arose from motivations in the field of geometric probability.   The authors are unaware of any similar research that considers the distribution of distances between points which themselves are distributed according to a Standard Normal Distribution in $k$-dimensional space without geometric boundaries.  This is surprising given the frequent application of the Standard Normal distribution for modeling purposes and the standard practice of normalizing data prior to the application of statistical learning techniques in order to prevent dominance of one variable over another.  In particular, the technique regarding normalization of data is highlighted in several of the key texts involving statistical learning. [8, 12]  This paper develops the distribution which defines the distances between randomly generated Standard Normal points in $k$-dimensional space without geometric boundaries, thus providing some statistical tools that support further refinement of machine learning and data mining applications.

\section{Distribution of Absolute Distance Between Two Gaussian Variables}

    This paper provides an investigation into the effects of dimensionality on the distribution of distances between random points which are distributed Standard Normal in each of their respective dimensions.  Most of the datasets confronted in high-dimensional problems are usually standardized in each of the attributes for respective datapoints. [8, 12]  Demonstrated below are the quantitative effects that dimensionality brings as the number of dimensions of the dataset increases.  To the authors' knowledge this subject has not been treated from a quantitative aspect in the literature. \\
\indent
    We begin by determining the absolute difference between two random variables ($p$ and $q$) which are both drawn from the standardized Gaussian distribution.  In this case, the absolute difference corresponds to one dimensional distance between points.  Note that this distance cannot be negative.  The Cumulative Density Function (CDF) which provides the distribution of the absolute difference between points $p$ and $q$, or their distance, is shown below and diagrammed further in Figure 1.

$$\int _{-\infty }^{\infty }\int _{p-x}^{p+x}\frac{e^{-\frac{p^2}{2}}
   e^{-\frac{q^2}{2}}}{\sqrt{2 \pi } \sqrt{2 \pi }}dqdp$$

\begin{center}
\includegraphics[scale=0.42]{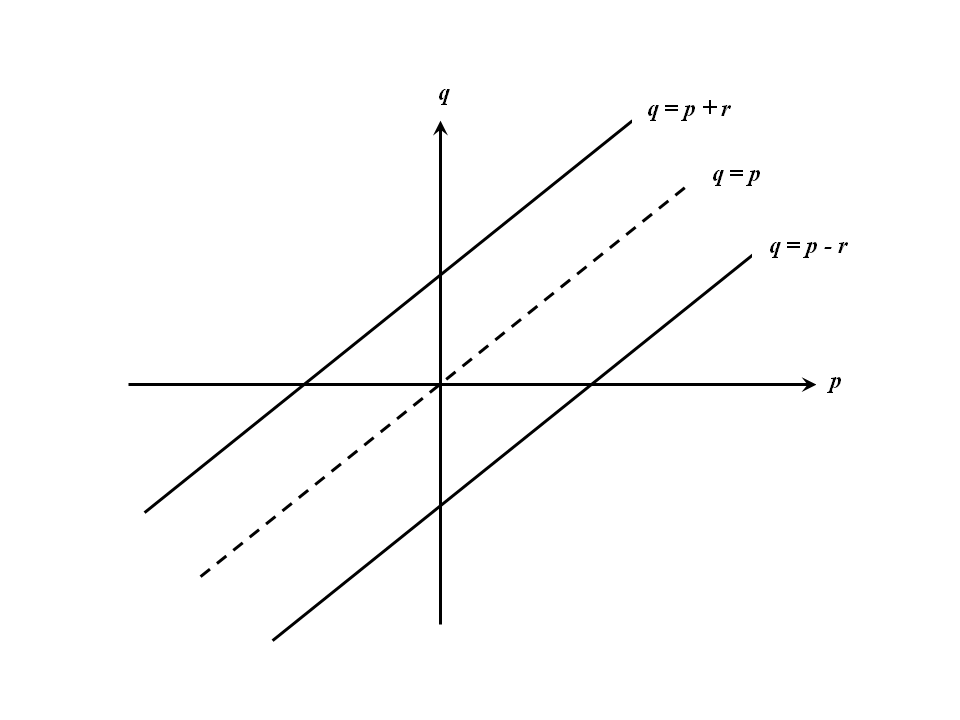}
\end{center}
Figure 1: Diagram of zone of integration which covers the probability that the distance between $p$ and $q$ is less than $r$

\indent
The inner term of the integral is the product of the Gaussian Probability Density Function (PDF) for each variable under consideration.  Integrating with respect to the innermost integral we arrive at the expression:
$$\int_{-\infty }^{\infty } \frac{e^{-\frac{p^2}{2}} \left(\int_0^{\frac{p+x}{\sqrt{2}}} e^{-t^2}
   \, dt-\int_0^{\frac{p-x}{\sqrt{2}}} e^{-t^2} \, dt\right)}{2 \sqrt{2 \pi }} \, dp$$
\indent
To determine the PDF of the above function, take the derivative of the integral with respect to $x$: %
\begin{align*}
f(x) &= \frac{\partial}{\partial x}\left[\int_{-\infty }^{\infty } \frac{e^{-\frac{p^2}{2}} \left(\int_0^{\frac{p+x}{\sqrt{2}}} e^{-t^2}
   \, dt-\int_0^{\frac{p-x}{\sqrt{2}}} e^{-t^2} \, dt\right)}{2 \sqrt{2 \pi }} \, dp\right] \\
&= \frac{e^{-\frac{x^2}{4}}}{\sqrt{\pi }}.
\end{align*}
\indent
Thus, the above expression is the PDF for the distribution of the absolute distance between two standardized Gaussian variables in one dimension.  A plot of this distribution is shown in Figure 2.
\begin{center}
\includegraphics[scale=0.3]{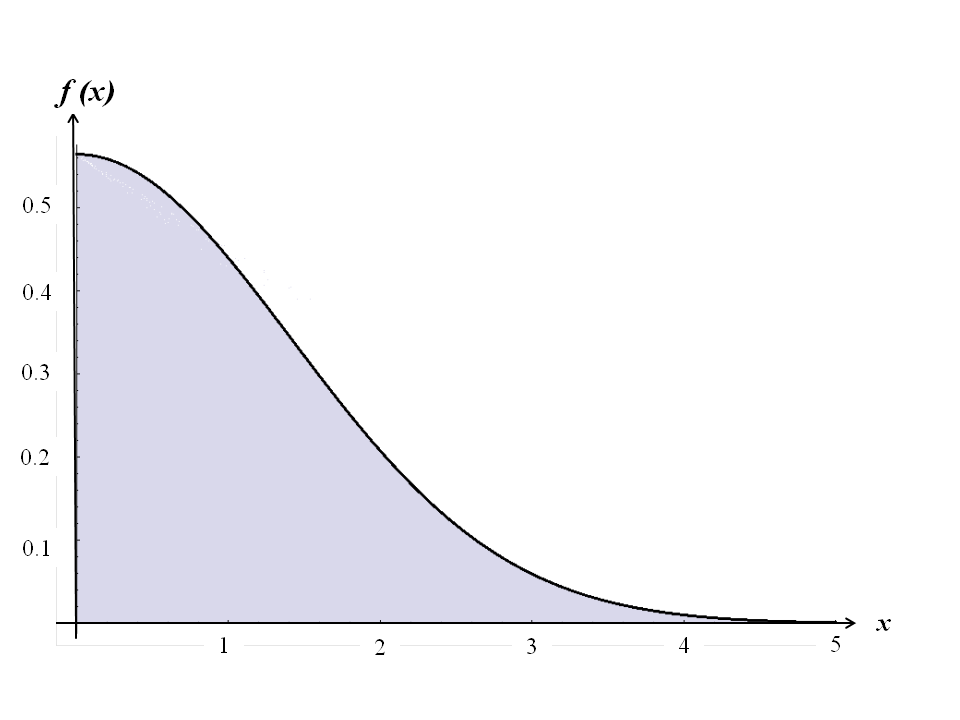} \\
Figure 2: Probability Density of the Absolute Difference between Points
\end{center}

\indent
In this section we have derived the PDF for the absolute difference (i.e. the non-negative distance) between any two random variables with a standard Gaussian distribution.  In the following sections, this concept is extended to an arbitrary number of $k$-dimensions.
\section{Extension of Distance to $k$ Dimensions of Gaussian Variables}
\indent The previous section derived the PDF for the distribution of the absolute difference in standard normal random variables in one dimension.  We will now extend the concept to an arbitrary $k$ number of dimensions and use the Euclidean metric to determine the distance between points.  Let two arbitrary points of interest in $k$-dimensions be further defined as $\Psi = \left(\psi_1 , \psi_2 , ...,\psi_k\right)$ and $\Gamma = \left(\gamma_1, \gamma_2, ...,\gamma_k\right)$.   We are interested in the Euclidean distance between the two points, which is defined as:

$$\left[\sum _{i=1}^k \left(\psi _i-\gamma _i\right){}^2\right]{}^{1/2}$$
\\ \\
\indent We generalize to $k$ dimensions now and begin by constructing the CDF which measures the probability that two points are separated by at most a distance less than or equal to $x$.  We do this by making a change of coordinates from standard Euclidan $x_i$ to those in a $k$-dimensional hypersphere and then concern ourselves only with the first $2^k$-tant where all of the values in the Euclidean coordinate system have positive values.  This is possible since each Euclidean coordinate is strictly non-negative and it has a corresponding probability that the difference between the two randomly distributed components is less than or equal to its radial value.  By approaching the problem in this manner, we need only concern ourselves with integration of all angles in the hypersphere from 0 to $\frac{\pi}{2}$, which we will designate using $\phi_i$ for each of the respective angular dimensions.  The radial component $R$ captures the probability that the CDF is less than or equal to $R$.  

\begin{center}
\includegraphics[scale=0.3]{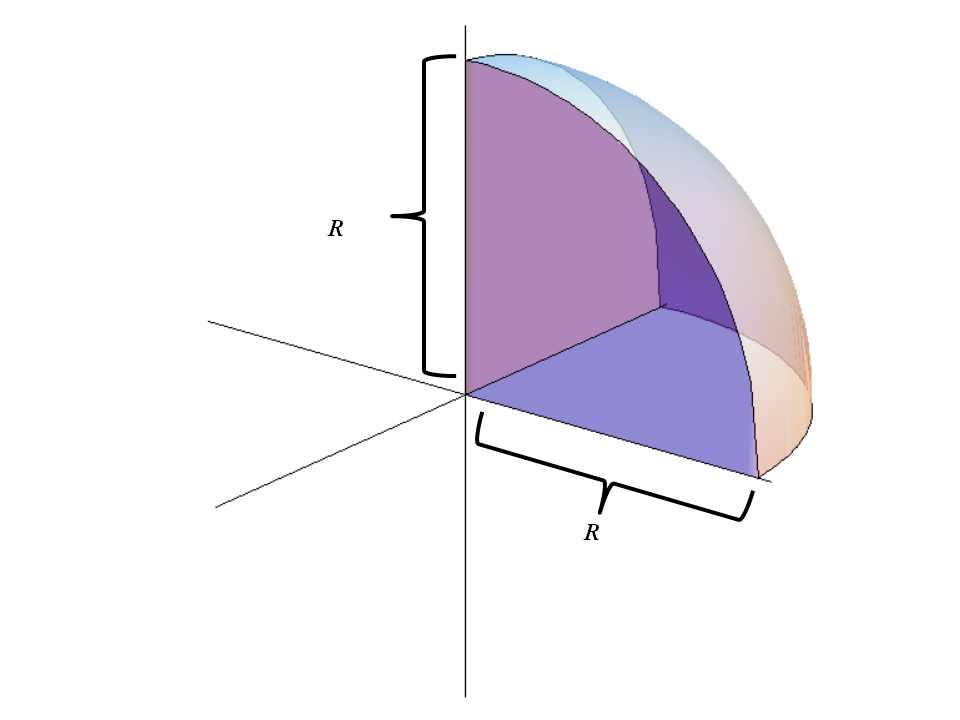}

Figure 3: Plot of region of iterated integration in the $2^k$-quadrant
\end{center}

We begin by making the following replacements of hyperspherical coordinates, as shown by Cohl in [5].
\begin{align*}
x_1 &= r \cos\left(\phi_1\right) \\
x_2 &= r \sin\left(\phi_1\right)\cos\left(\phi_2\right) \\
x_3 &= r \sin\left(\phi_1\right)\sin\left(\phi_2\right)\cos\left(\phi_3\right) \\
\vdots \\
x_{k-1} &=  r \sin\left(\phi_1\right)...\sin\left(\phi_{k-2}\right)\cos\left(\phi_{k-1}\right) \\
x_k &=  r \sin\left(\phi_1\right)...\sin\left(\phi_{k-2}\right)\sin\left(\phi_{k-1}\right) 
\end{align*}
After converting to hyperspherical coordinates and integrating, the spherical volume element $d^k V$ is given by
$$d^k V = r^{k-1} \sin^{k-2}\left(\phi_1\right)\sin^{k-3}\left(\phi_2\right)...\sin\left(\phi_{k-2}\right) dr \  d\phi_1 \ d\phi_2 ... d\phi_{k-1}.$$
The iterated integral in the first $2^k$- tant which defines the CDF for the probability that the distance between two random Gaussian points is no greater than $R$ in $k$ dimensional space is then given by
$$F(R,k) = \int^{\pi/2}_0...\int^{\pi/2}_0\int^R_0 ( \frac{e^{-\frac{1}{4} (r \cos \left( \phi _1 \right) ){}^2}}{\sqrt{\pi}} ) ( \frac{e^{-\frac{1}{4} (r \sin \left( \phi _1 \right) \cos \left( \phi_2 \right) ){}^2}}{\sqrt{\pi}} )
( \frac{e^{-\frac{1}{4} (r \sin \left( \phi _1 \right) \sin \left( \phi_2 \right) \cos \left( \phi_3 \right) ){}^2}}{\sqrt{\pi}} ) ...$$ \\ 
$$...( \frac{e^{-\frac{1}{4} (r \sin \left( \phi _1 \right)...\sin \left( \phi_{k-2} \right) \cos \left( \phi_{k-1} \right) ){}^2}}{\sqrt{\pi}} )
( \frac{e^{-\frac{1}{4} (r \sin \left( \phi _1 \right)...\sin \left( \phi_{k-2} \right) \sin \left( \phi_{k-1} \right) ){}^2}}{\sqrt{\pi}} )... $$ \\
$$...r^{k-1} \sin^{k-2}\left(\phi_1\right)\sin^{k-3}\left(\phi_2\right)...\sin\left(\phi_{k-2}\right) dr \ d\phi_1 \ d\phi_2...d\phi_{k-1}.
 $$  \\
\indent However, since we have completed the change of coordinates to the first $2^k$-tant of a hypersphere, the above expansion simplifies to \\
$$F(R,k) = \int^{\pi/2}_0...\int^{\pi/2}_0\int^R_0 \left(\frac{1}{\sqrt{\pi}}\right)^k\left(e^{-\frac{r^2}{4}}\right)
r^{k-1} \sin^{k-2}\left(\phi_1\right)\sin^{k-3}\left(\phi_2\right)...$$
$$...\sin\left(\phi_{k-2}\right) dr \ d\phi_1 \ d\phi_2...d\phi_{k-1}.
$$
\indent Rearranging the terms so that we are only integrating with respect to the variable in question at each iterated integral, we have
$$F(R, k) = \left(\frac{1}{\sqrt{\pi}}\right)^k \int^{\pi/2}_0(\int^{\pi/2}_0 (\sin(\phi_{k-2}))\dots\int^{\pi/2}_0 ( (\sin^{k-3}(\phi_2)) \dots $$
$$\int^{\pi/2}_0 ( (\sin^{k-2}(\phi_1))(\int^R_0(e^{-\frac{r^2}{4}}) r^{k-1} dr)...d\phi_1 ) \  d\phi_2 ) ...\\
d\phi_{k-2} ) \  d\phi_{k-1}.
$$
\indent    Beginning with the innermost integral, we make use of the solution to the integral
$$\int^R_0 e^{-\frac{r^2}{4}} r^{k-1}dr = 2^{k-1}R^k\left(R^2 \right)^{-k/2}
\left( \Gamma\left( \frac{k}{2}\right)-\Gamma\left( \frac{k}{2},\frac{R^2}{4}\right)\right).$$
\indent  Since $k$ is a positive integer, the above result simplifies to   \\
$$\int^R_0 e^{-\frac{r^2}{4}} r^{k-1}dr = 2^{k-1}\left( \Gamma\left( \frac{k}{2}\right)-\Gamma\left( \frac{k}{2},\frac{R^2}{4}\right)\right). $$
\indent    After integration, the result is a constant with respect to all subsequent iterated integrals.  Thus, we can move it to the outside of the iterated integrals as shown below.
$$F(R,k) = \left( \frac{1}{\sqrt{\pi}}\right)^k \left( 2^{k-1}\left( \Gamma\left( \frac{k}{2}\right)-\Gamma\left( \frac{k}{2},\frac{R^2}{4}\right)\right)\right)\dots $$
$$...\int^{\pi/2}_0(\int^{\pi/2}_0 (\sin(\phi_{k-2}))... \\ \\
\int^{\pi/2}_0 ( (\sin^{k-3}(\phi_2))\int^{\pi/2}_0 ( \sin^{k-2}(\phi_1)
 d\phi_1 ) \  d\phi_2 )... $$
$$...d\phi_{k-2} ) \  d\phi_{k-1}.$$ \\
\indent   Now we have a series of iterated integrals for decreasing powers of the sine function for all integral values for $k-2$ to $1$.  We now use the following identity which holds for any $m$ such that the real component is greater than $-1$ and can be verified using a computer algebra system, 
$$\int^{\pi/2}_0\sin^m(\theta) d\theta = \frac{\sqrt{\pi} \ \Gamma\left(\frac{m+1}{2}\right)}
{2 \ \Gamma\left(\frac{m}{2}+1\right)}.$$

Since all powers of $m$ in the integral are strictly positive integers, each integral in our iterated integral expression can be replaced with the closed form expression above.  The result is strictly a constant depending only on the power that the sine function is raised to in the integrand.  Since integration occurs for every power of sine for $1 \le m \le k-2$, the results of each of these integrals can be gathered as constants to the outside of the integral and formed into a product extending from $1$ to $k-2$ which leaves us with the following result
$$F(R,k) = \left( \frac{1}{\sqrt{\pi}}\right)^k \left( 2^{k-1}\left( \Gamma\left( \frac{k}{2}\right)-\Gamma\left( \frac{k}{2},\frac{R^2}{4}\right)\right)\right) 
\left( \prod^{k-2}_{m=1}\ \frac{\sqrt{\pi}\Gamma\left(\frac{m+1}{2}\right)}
{2\Gamma\left(\frac{m}{2}+1\right)}\right)
\int^{\pi/2}_0d\phi_{k-1}.$$
\indent   The final integral found on the far right evaluates to $\frac{\pi}{2}$, which leaves us the simplified CDF shown below. 
$$F(R,k) = \left( \frac{1}{\sqrt{\pi}}\right)^k \left( 2^{k-1}\left( \Gamma\left( \frac{k}{2}\right)-\Gamma\left( \frac{k}{2},\frac{R^2}{4}\right)\right)\right) 
\left( \prod^{k-2}_{m=1}\ \frac{\sqrt{\pi}\Gamma\left(\frac{m+1}{2}\right)}
{2\Gamma\left(\frac{m}{2}+1\right)}\right)
\left( \frac{\pi}{2}\right).$$ 
\indent After further simplification this reduces to
$$F(R,k)=\left( \Gamma\left( \frac{k}{2}\right)-\Gamma\left( \frac{k}{2},\frac{R^2}{4}\right)\right)
\prod^{k-2}_{m=1}\frac{\Gamma\left( \frac{m+1}{2}\right)}{\Gamma\left(\frac{m}{2}+1\right)}.$$
\indent Closer examination of the iterated product by expansion of terms and cancellation produces the following identity:
$$\prod^{k-2}_{m=1}\frac{\Gamma\left( \frac{m+1}{2}\right)}{\Gamma\left(\frac{m}{2}+1\right)} = \frac{\Gamma\left(1\right)}{\Gamma\left(\frac{k}{2}\right)} =
\frac{1}{\Gamma\left(\frac{k}{2}\right)}.$$
\indent
Substituting this result back into the CDF, we have
\begin{align*}
F(R,k)&=\left( \Gamma\left( \frac{k}{2}\right)-\Gamma\left( \frac{k}{2},\frac{R^2}{4}\right)\right)
\prod^{k-2}_{m=1}\frac{\Gamma\left( \frac{m+1}{2}\right)}{\Gamma\left(\frac{m}{2}+1\right)} \\ 
&= \frac{\left(\Gamma\left(\frac{k}{2}\right)-\Gamma\left(\frac{k}{2},\frac{R^2}{4}\right)\right)}{\Gamma\left(\frac{k}{2}\right)} \\
&=1-\frac{\Gamma\left(\frac{k}{2},\frac{R^2}{4}\right)}{\Gamma\left(\frac{k}{2}\right)}.
\end{align*}
\indent Taking the derivative with respect to $R$ yields the PDF below.
\begin{align*}
f(R,k)&=\frac{\partial}{\partial R}\left[1-\frac{\Gamma\left(\frac{k}{2},\frac{R^2}{4}\right)}{\Gamma\left(\frac{k}{2}\right)}\right] \\ 
&= \frac{2^{1-k}e^{-\frac{R^2}{4}}R^{k-1}}{\Gamma\left(\frac{k}{2}\right)}.
\end{align*}
\indent This reduces to $f(R, k)  =  \frac{2^{1-k}e^{-\frac{R^2}{4}}R^{k-1}}{\Gamma\left(\frac{k}{2}\right)}$, where $R$ is the absolute difference between two random variables which have each of their $k$ components distributed Gaussian.  This distribution is strictly non-negative, and a quick integration over all possible values of $R$ from 0 to infinity equals 1, thereby confirming the validity of the PDF. 

\indent A plot of the PDF is shown in Figure 4 for a variety of dimensions ($k$), where $k=1, 2, 3, 4, 5, 10, 20, 30, 40, 50, 100$.  This plot is interesting in that as $k \rightarrow \infty$, the distribution resembles a shifted normal.  Furthermore, due to the underlying calculation of Euclidean distance, there is no asymptotic limit regarding shift occurence. \\

\begin{center}
\includegraphics[scale=0.6]{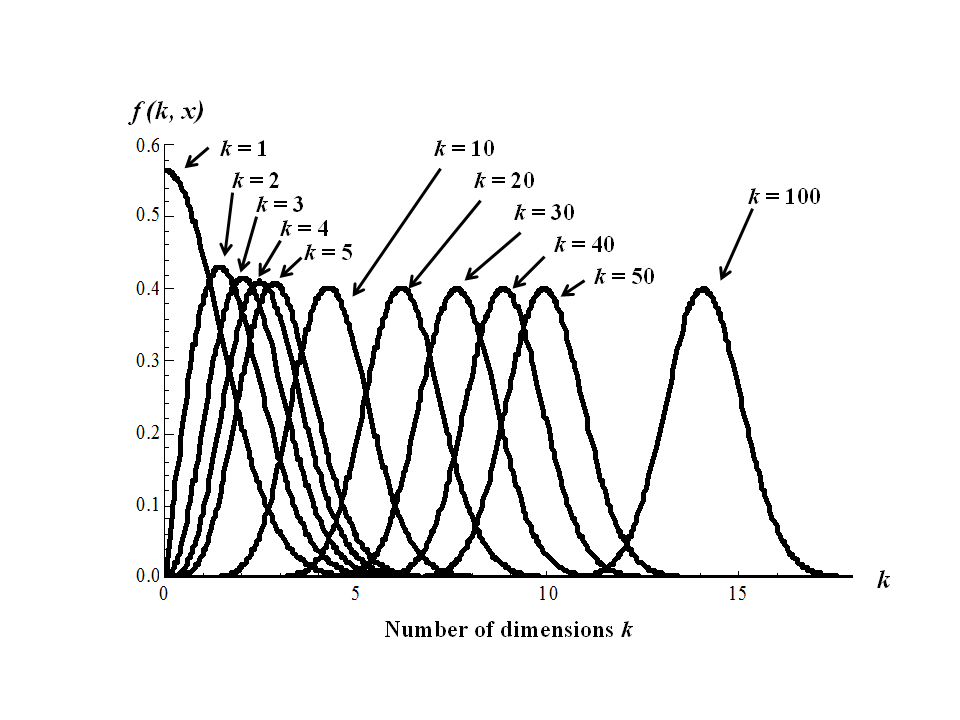}
Figure 4: Series of overlaid plots of PDFs for various values of $k$
\end{center}

\indent    In summary, we evaluated the Euclidean distance between the absolute differences of two $k$-dimensional random variables by examining them as the equivalent of a set of points in the first $2^n$-tant of a hypersphere.  Through this conversion, we are able to evaluate an iterated integral to a simple closed form expression for the CDF which reflects the distance between two such points.  Further differentiation with respect to the radial component $r$ produces the PDF of the distribution for the distance between two points of interest in $k$-dimensional space.

\newpage

\section{Calculation of Raw and Central Moments}

\indent Anytime that a new distribution is derived, it is natural to discuss the moments.  We will do this by evaluating the following integral to calculate the $n^{th}$ raw moment of the distribution. \\
$$\int^\infty_0 R^n\left(\frac{2^{1-k}e^{-\frac{R^2}{4}}R^{k-1}}{\Gamma\left(\frac{k}{2}\right)}\right) dR$$
\indent
Factoring out constants, we are left with
$$\frac{2^{1-k}}{\Gamma\left(\frac{k}{2}\right)}
\int^\infty_0\left(e^{-\frac{R^2}{4}}R^{n+k-1}\right)dR.$$
\indent
The following identity holds, provided that $k+n$ is strictly greater than zero.  This is always true in this analysis given that the dimensionality has $k\ge1$ and the raw moment has $n\ge1$.
$$\int^\infty_0\left(e^{-\frac{R^2}{4}}R^{n+k-1}\right)dR=2^{n+k-1}\Gamma\left(\frac{k+n}{2}\right).$$
\indent This simplifies to the expression below for the $n^{th}$ raw moment in $k$-dimensions.
$$\int^\infty_0 R^n\left( \frac{2^{1-k}e^{-\frac{R^2}{4}}R^{k-1}}{\Gamma\left(\frac{k}{2}\right)}\right) dR
=\frac{2^n\Gamma\left(\frac{k+n}{2}\right)}{\Gamma\left(\frac{k}{2}\right)}.$$
\indent Therefore, the first four raw moments for $k$-dimensions are shown below. \\
\\
\indent
$\int^\infty_0 R\left( \frac{2^{1-k}e^{-\frac{R^2}{4}}R^{k-1}}{\Gamma\left(\frac{k}{2}\right)}\right) dR
=\frac{2\Gamma\left(\frac{k+1}{2}\right)}{\Gamma\left(\frac{k}{2}\right)}$  (first raw moment). \\ \\
\indent $\int^\infty_0 R^2\left( \frac{2^{1-k}e^{-\frac{R^2}{4}}R^{k-1}}{\Gamma\left(\frac{k}{2}\right)}\right) dR
=\frac{4\Gamma\left(\frac{k+2}{2}\right)}{\Gamma\left(\frac{k}{2}\right)}$  (second raw moment). \\ \\
\indent $\int^\infty_0 R^3\left( \frac{2^{1-k}e^{-\frac{R^2}{4}}R^{k-1}}{\Gamma\left(\frac{k}{2}\right)}\right) dR
=\frac{8\Gamma\left(\frac{k+3}{2}\right)}{\Gamma\left(\frac{k}{2}\right)}$  (third raw moment). \\ \\
\indent $\int^\infty_0 R^4\left( \frac{2^{1-k}e^{-\frac{R^2}{4}}R^{k-1}}{\Gamma\left(\frac{k}{2}\right)}\right) dR
=\frac{16\Gamma\left(\frac{k+4}{2}\right)}{\Gamma\left(\frac{k}{2}\right)}$  (fourth raw moment).

\indent Calculation of the central moments begins with the consideration of the first raw moment identified above and is used to calculate the $n^{th}$ central moment in the integral
$$\int^\infty_0\left( x-\mu\right)^n f(x) dx.$$
\indent Substituting values for $\mu$ and $f(x)$ and restricting the lower bound of the integral to reflect strictly non-negative values, we have the following expression.
$$\int^\infty_0\left(R-\frac{2\Gamma\left(\frac{k+1}{2}\right)}{\Gamma\left(\frac{k}{2}\right)}\right)^n \left(\frac{2^{1-k}e^{-\frac{R^2}{4}}R^{k-1}}{\Gamma\left(\frac{k}{2}\right)}\right) dR$$
\indent Unfortunately, determining a closed form expression for the $n^{th}$ central moment is non-trivial, but this paper will consider the second, third, and fourth central moments.  We begin by determining the second central moment $\mu_2$. \\
$$\mu_2 = \int^\infty_0\left(R-\frac{2\Gamma\left(\frac{k+1}{2}\right)}{\Gamma\left(\frac{k}{2}\right)}\right)^2 \left(\frac{2^{1-k}e^{-\frac{R^2}{4}}R^{k-1}}{\Gamma\left(\frac{k}{2}\right)}\right) dR = 2k-\frac{4\Gamma\left(\frac{k+1}{2}\right)^2}{\Gamma\left(\frac{k}{2}\right)^2} $$
\\
\indent Now determine the third central moment $\mu_3$:
\begin{align*}\mu_3 &= \int^\infty_0\left(R-\frac{2\Gamma\left(\frac{k+1}{2}\right)}{\Gamma\left(\frac{k}{2}\right)}\right)^3 \left(\frac{2^{1-k}e^{-\frac{R^2}{4}}R^{k-1}}{\Gamma\left(\frac{k}{2}\right)}\right) dR \\
&= 
\frac{16\Gamma\left(\frac{k+1}{2}\right)^3 + 4(1-2k)\Gamma\left(\frac{k}{2}\right)^2\Gamma\left(\frac{k+1}{2}\right)}{\Gamma\left(\frac{k}{2}\right)^3}
\end{align*}
\indent Finally, determine the fourth central moment $\mu_4$:
\begin{align*}\mu_4 &= \int^\infty_0\left(R-\frac{2\Gamma\left(\frac{k+1}{2}\right)}{\Gamma\left(\frac{k}{2}\right)}\right)^3 \left(\frac{2^{1-k}e^{-\frac{R^2}{4}}R^{k-1}}{\Gamma\left(\frac{k}{2}\right)}\right) dR \\
&= 4k(k+2)+
\frac{\pi4^{3-k}(k-2)\Gamma\left(k\right)^2-48\Gamma\left(\frac{k+1}{2}\right)^4}
{\Gamma\left(\frac{k}{2}\right)^4} 
\end{align*}
\\
\indent The skewness for this distribution is as indicated below using $\mu_i$, where $i$ is the $i$-th central moment.
$$\gamma_1 = \frac{\mu_3}{\mu^{3/2}_2} = \frac{\sqrt{\pi} \ 2^{\frac{3}{2}-k } \ \Gamma(k)((1-2k) \ \Gamma(\frac{k}{2})^2+4 \ \Gamma(\frac{k+1}{2})^2)}{\Gamma(\frac{k}{2}) \ (k \ \Gamma(\frac{k}{2})^2 - 2 \ \Gamma(\frac{k+1}{2})^2)^\frac{3}{2}}$$
\\
\indent The kurtosis for this distribution is provided below using the expression for $\beta_2$.
$$\beta_2 = \frac{\mu_4}{\mu^{2}_2} = 4k(k + 2) + \frac{\pi \ 4^{3-k}(k-2) \  \Gamma(k)^2 - 48 \ \Gamma(\frac{k+1}{2})^4}{\Gamma(\frac{k}{2})^4}
$$
\\
\indent All of the above identities hold provided that $Re(k) \ge 0$, which holds for the problem under consideration in this paper.  Thus, the closed form expression for the $n^{th}$ raw moment of the distance distribution in $k$-dimensions is shown.  While a closed form expression for the $n^{th}$ central moment is not derived, the second, third, and fourth central moments are directly calculated. \\

\section{Discussion and Further Work}
In this paper we have developed the PDF for the distribution of Euclidean distance between any two points with $k$ features or dimensions.  Our methodology utilized a change of coordinates and allowed for the derivation of the CDF for the distribution of distances between points in $k$-dimensions.  Subsequent derivation provided the pdf and we were pleasantly surprised with the elegant closed-form expression.

The distribution of distances between random points is especially important to data mining applications since several important classification algorithms rely on the use of distance to nearest known examples to determine the classification of unknown instances.  Frequently it is the case that a significant portion of features or dimensions in a given dataset are completely irrelevant to the problem, or are themselves noise. The inclusion of variables which provide no additional accuracy actually serves to increase the distance between points, and the effect of the added distance is distributed as outlined in this paper, and is dependent upon the number of extraneous dimensions.  Additionally, spurious results may be given erroneous meaning when such results emerge from the inclusion of a large number of randomly distributed features in the data under consideration.  

Since we have derived the distribution of distances between neighboring points, future research will evaluate the use of this distribution in the reduction of such spurious associations and the probability that these associations would arise from the number of underlying features.

\section{References}

\noindent [1] Aggarwal, Charu C., Alexander Hinneburg, and Daniel A. Keim. On the surprising behavior of distance metrics in high dimensional space. Springer Berlin Heidelberg, 2001.

\noindent [2] Alagar, Vangalur S. ``The distribution of the distance between random points." Journal of Applied Probability (1976) : 558 - 566.

\noindent [3] Bellman, Richard Ernest. Rand Corporation (1957). Dynamic programming. Princeton University Press. ISBN 978-0-691-07951-6.,
Republished: Richard Ernest Bellman (2003). Dynamic Programming. Courier Dover Publications. ISBN 978-0-486-42809-3.

\noindent [4] Burgstaller, Bernhard, and Friedrich Pillichshammer. ``The average distance between two points." Bulletin of the Australian Mathematical Society 80.03 (2009): 353-359.

\noindent [5] Cohl, Howard S. ``Fundamental Solution of Laplace's Equation in Hyperspherical Geometry." SIGMA. Symmetry, Integrability and Geometry: Methods and Applications 7 (2011): 108.

\noindent [6] Hammersley, John M. ``The distribution of distance in a hypersphere." The Annals of Mathematical Statistics (1950): 447-452.

\noindent [7] Hastie, Trevor, et al. ``The elements of statistical learning." Vol. 2. No. 1. New York: Springer, 2009.

\noindent [8] Larose, Daniel T. ``Discovering knowledge in data: an introduction to data mining. 2005." (2003).

\noindent [9] Lord, Reginald Douglas. ``The distribution of distance in a hypersphere." The Annals of Mathematical Statistics (1954): 794-798.

\noindent [10] Saff, Edward B., and A. BJ Kuijlaars. ``Distributing many points on a sphere." The Mathematical Intelligencer 19.1 (1997): 5-11.

\noindent [11] Solomon, Herbert. Geometric probability. Vol. 28. SIAM, 1978.

\noindent [12] Tan, Pang-Ning, Michael Steinbach, and Vipin Kumar. ``Introduction to Data Mining." (2006).

\end{document}